\documentclass[12pt]{article}

\usepackage[T1]{fontenc}
\usepackage[latin1]{inputenc}
\usepackage{bbm}
\usepackage{stmaryrd}
\usepackage{latexsym}
\usepackage{amsfonts}
\usepackage{amsmath,amssymb,amscd}
\usepackage{mathtools}
\usepackage{yfonts}
\usepackage{dsfont}
\usepackage{mathabx}
\usepackage[dvips]{graphicx}
\usepackage{epsfig}
\usepackage{psfrag}
\usepackage[font={small}]{caption}
\usepackage{color}
\usepackage{float}
\usepackage{upgreek}
\usepackage{tikz}
\usepackage{tikz-cd}
\usepackage{relsize}
\usepackage{url}
\DeclareFontFamily{U}{txsyc}{}
\DeclareFontShape{U}{txsyc}{m}{n}{
   <-> txsyc%
}{}
\DeclareFontShape{U}{txsyc}{bx}{n}{
   <-> txbsyc%
}{}
\DeclareFontShape{U}{txsyc}{l}{n}{<->ssub * txsyc/m/n}{}
\DeclareFontShape{U}{txsyc}{b}{n}{<->ssub * txsyc/bx/n}{}
\DeclareSymbolFont{symbolsC}{U}{txsyc}{m}{n}
\SetSymbolFont{symbolsC}{bold}{U}{txsyc}{bx}{n}
\DeclareFontSubstitution{U}{txsyc}{m}{n}
\DeclareMathSymbol{\df}{\mathrel}{symbolsC}{"42}
\DeclareMathSymbol{\fd}{\mathrel}{symbolsC}{"43}
\DeclareMathSymbol{\lJoin}{\mathrel}{symbolsC}{"58}
\DeclareMathSymbol{\rJoin}{\mathrel}{symbolsC}{"59}

\newcommand{\f}[2]{\frac{#1}{#2}}
\newcommand{\cA}{\mathcal{A}}

\newcommand{\cC}{\mathcal{C}}

\newcommand{\cE}{\mathcal{E}}

\newcommand{\cL}{\mathcal{L}}

\newcommand{\cM}{\mathcal{M}}

\newcommand{\cU}{\mathcal{U}}

\newcommand{\EE}{\mathbb{E}}

\newcommand{\NN}{\mathbb{N}}
\newcommand{\PP}{\mathbb{P}}

\newcommand{\RR}{\mathbb{R}}
\renewcommand{\SS}{\mathbb{S}}

\newcommand{\ff}{\mathfrak{f}}

\newcommand{\iy}{\infty}

\newcommand{\lt}{\left}
\newcommand{\me}{\medskip}
\newcommand{\na}{\nabla}
\newcommand{\pa}{\partial}
\newcommand{\ri}{\rightarrow}
\newcommand{\rt}{\right}
\newcommand{\sm}{\smallskip}
\newcommand{\tr}{\triangle}

\newcommand{\wi}{\widetilde}

\newcommand{\dive}{\mathrm{div}}

\newcommand{\fo}{\forall\ }

\newcommand{\Hess}{\mathrm{Hess}\,}

\newcommand{\lan}{\lt\langle}

\newcommand{\ran}{\rt\rangle}

\newcommand{\st}{\,:\,}

\newcommand{\trace}{\mathrm{tr}}

\newcommand{\bq}{\begin{eqnarray*}}
\newcommand{\bqn}[1]{\begin{eqnarray}\label{#1}}
\newcommand{\eq}{\end{eqnarray*}}
\newcommand{\eqn}{\end{eqnarray}}
\newcommand{\wwtbp}{\hfill $\blacksquare$\par\me\noindent}

\newcommand{\thistitlepagestyle}{}
\newcommand{\lin}{\llbracket}
\newcommand{\rin}{\rrbracket}

\newcommand{\ttsim}{\raise.17ex\hbox{$\scriptstyle\mathtt{\sim}$}}

\newcommand{\kh}{\kern .08em}

\newtheorem{pro}{Proposition}

\newtheorem{lem}[pro]{Lemma}
\newtheorem{theo}[pro]{Theorem}

\renewcommand{\thepro}{\arabic{pro}}

\newenvironment{rem}
{\par\me\refstepcounter{pro}\noindent{\bf Remark \thepro\ }}
{\hfill $\square$\par\sm\noindent}

\newcommand{\proof}{\par\me\noindent\textbf{Proof}\par\sm\noindent}

\newcommand{\sU}{\textsf{U}}
\newcommand{\sV}{\textsf{V}}
\newcommand{\sL}{\textsf{L}}
\newcommand{\sX}{\textsf{X}}
\newcommand{\sH}{\textsf{H}}
 \newcommand{\eb}{^{(\beta)}}
 \newcommand{\ebx}{^{(\beta,x)}}

\title{The asymptotic behavior of fraudulent algorithms}
\author{Michel Benaïm${}^\dagger$ and Laurent Miclo${}^\ddagger$
}
 \date{\vbox{\copy1
 \vskip5mm
 \copy2
}
}

\begin{document}

\setbox1=\vbox{
 \large
 \begin{center}
 ${}^\dagger$ 
 Institut de Mathématiques \\
 Neuchâtel University
 \end{center}
 }

\setbox2=\vbox{
\large
\begin{center}
 ${}^\ddagger$
Toulouse School of Economics\\
Institut de Mathématiques de Toulouse\\
CNRS and University of Toulouse
\end{center}
}
\setbox5=\vbox{
\hbox{miclo@math.cnrs.fr\\[1mm]}
\vskip1mm
\hbox{Toulouse School of Economics\\}
\hbox{1, Esplanade de l'université\\}
\hbox{31080 Toulouse cedex 6, France\\[1mm]}
\hbox{Institut de Mathématiques de Toulouse\\}
\hbox{Université Paul Sabatier, 118, route de Narbonne\\}
\hbox{31062 Toulouse cedex 9, France.\\[1mm]}
}
\setbox4=\vbox{
\hbox{michel.benaim@unine.ch\\[1mm]}
\vskip1mm
\hbox{Institut de Mathématiques}
\hbox{Faculté des sciences, Université de Neuchâtel}
\hbox{11 rue E. Argand, 2000 Neuchâtel, Suisse}
}

 \maketitle
\thistitlepagestyle
\abstract{
Let $U$ be a Morse function on a  compact connected $m$-dimensional Riemannian manifold, $m \geq 2,$ satisfying $\min U=0$ and let $\cU = \{x \in M \: : U(x) = 0\}$ be the set of global minimizers. Consider the stochastic algorithm $X\eb\df(X\eb(t))_{t\geq 0}$ defined on $N = M \setminus \cU,$  whose generator is
$U \tr \cdot-\beta\lan \na U,\na \cdot\ran$, where $\beta\in\RR$ is a real parameter.
We show that for $\beta>\frac{m}{2}-1,$ $X\eb(t)$ converges a.s.\ as $t \rightarrow   \infty,$ toward a point $p \in \cU$ and that each  $p \in \cU$  has a positive probability to be selected. On the other hand, for $\beta <  \frac{m}{2}-1,$ the law of $(X\eb(t))$ converges in total variation (at an exponential rate) toward the probability measure  $\pi_{\beta}$ having density proportional to  $U(x)^{-1-\beta}$ with respect to the Riemannian measure.
}
\vfill\null
{\small
 \par\vskip.3cm
 \textbf{Fundings:} the grants SNF 200020-219913, ANR-17-EURE-0010 and AFOSR-22IOE016  are acknowledged.
}\par

\newpage

\section{Introduction}

Stochastic global minimization algorithms taking into account the a priori  knowledge of the minimal value of the objective function $U$  are called \textbf{fraudulent}, since this minimal value is often  not available in practice.
Nevertheless some of their interests are presented in \cite{miclo:hal-04094950},
where such a procedure was introduced  when $U$ is a Morse function defined on a compact manifold $M$ of dimension $m\geq 2$.
The underlying stochastic process $X^{(\beta)}\df (X^{(\beta)}(t))_{t\geq 0}$, taking values in $M$, comes with a real parameter $\beta$ which can be tuned to increase the relative importance of $U$ with respect to the injected randomness.
Two quantities $\beta_\vee\geq \beta_\wedge \in\RR$ (depending explicitly on the eigenvalues of the Hessians of $U$ at its global minima, see Remark \ref{rem3} below) were introduced
so that $\beta>\beta_\vee$ implies the a.s.\ convergence of $X^{(\beta)}(t)$  as $t \rightarrow \infty,$  toward the global minima of $U$ (and each of them
attracts the algorithm with positive probability, when
$X^{(\beta)}(0)$ is not a global minima),
while for $\beta< \beta_\wedge$ the probability that  $X^{(\beta)}(t)$ a.s.\ converges  toward a global minimum of $U$ is zero.

Our goal here is to sharpen these result and describe completely the long term behavior of $X^{(\beta)}$ for all $\beta \neq \beta_0,$ where $\beta_0 \df \frac{m}{2}-1$ is a universal (i.e independent of $U$) critical value.  We will show that for $\beta > \beta_0,$ $X^{(\beta)}(t)$ a.s.\ converges  toward a global minimizer of $U$ and that each global minimizer has a positive probability to be selected. On the other hand, for $\beta < \beta_0,$ the process converges in distribution toward a (unique) invariant distribution whose density (with respect to the Riemannian measure) is explicit.
This result will be a consequence of the persistence/non-persistence approach presented in \cite{MR3910006} and \cite{Ben18}. The paper is organized as follows. Section \ref{sec:main} sets the notation and presents the main results. Section \ref{sec:Euclide} considers the situation where $M$ is no longer a compact manifold but the Euclidean space $\RR^m.$ It allows to introduce the main ingredients of the proof in a simple setting. Section \ref{sec:proofs} is devoted to the proof of the main results. Certain additional points are discussed in appendix.
\section{Notation and main result}
\label{sec:main}
 We assume throughout that $M$ is a compact connected Riemannian manifold having dimension $m \geq 2$ and $U : M \rightarrow \RR$ is a smooth function such that (this is the fraudulent assumption):
\bq
\min_M U&=&0.\eq
The zero set of $U,$ \bq
\cU&\df&\{p\in M\st U(p)=0\},\eq is then the set of global minimizers. We furthermore assume that every $p \in \cU$ is {\em non-degenerate}, meaning that the Hessian of $U$ at $p$ is non-degenerate.
This assumption implies that $\cU$ is finite. In particular, $N \df M \setminus \cU$ is  a noncompact connected manifold.
\par
Let $L_{\beta}$ be the operator on $\cC^2(M)$ defined as
\bqn{Lb}
L_\beta&\df&U \tr \cdot-\beta\lan \na U,\na \cdot\ran\eqn
where $\tr$, $\lan\cdot,\cdot\ran$ and $\na$ stand for the Laplacian, scalar product and gradient associated to the Riemannian structure of $M,$ and $\beta \in \RR.$

A {\em diffusion process generated by (\ref{Lb})}, is a continuous-time Feller Markov process on $M,$ $X\eb = (X\eb(t))_{t \geq 0},$ with infinitesimal generator $\cL_{\beta}$ and domain $D(\cL_{\beta}) \subset \cC^0(M)$ (see e.g.\ Le Gall \cite{MR3497465}, Section~6.2, for the definitions of Feller processes, domains and generators) such that for all  $f \in \cC^2(M):$

 $$f \in D(\cL_{\beta}) \mbox{ and } \cL_{\beta} f = L_{\beta} f.$$
Since the mapping $\na U$  and  $\sqrt{U}$ are Lipschitzian, due to
the non-degeneracy  assumption of the zeroes of $U$ for the latter, such a diffusion process exists. More details are given in the appendix. In addition, given the initial distribution of $X\eb(0),$  say $\mu,$   the law of $X\eb, \PP\eb_{\mu},$ is uniquely determined by $\mu$ and $L_{\beta}.$   As usual, we write  $\PP_x\eb$ for $\PP\eb_{\delta_x}.$ By a mild (but convenient) abuse of notation we may write $\PP_x(X\eb \in \cdot)$ for $\PP_x\eb(\cdot).$
We also let  $P\eb = (P_t\eb)_{t \geq 0}$ denote the semi-group induced by $X\eb.$ It is defined, as usual, by
$$\fo t\geq 0,\, \fo x\in M,\qquad P_t\eb f(x) \df \EE_x f(X\eb(t))$$ for every measurable, bounded or nonnegative, map $f : M \rightarrow \RR.$

The next proposition summarizes some basic properties of $P\eb.$ Its proof, which relies on classical results, is  given in the appendix.
 \begin{pro}\label{prosemigroup}
\begin{description}
\item[(i)]  $P\eb$ leaves $N$ and $\cU$ invariant:

- For all $t \geq 0,\, P_t\eb \mathbf{1}_N = \mathbf{1}_N.$
\item[(ii)] $P\eb$ is Feller on $M$ and  strong-Feller on $N:$

-  For all $t \geq 0$ and $f \in \cC^0(M),\, P\eb_t (f) \in  \cC^0(M);$

-  For all $t > 0$ and $f : N \rightarrow \RR$ bounded measurable,  $P_t\eb(f)$ is continuous on $N$.
\end{description}
\end{pro}
Note that $P\eb$ is not strong Feller on $M$, as it can be seen by considering the indicator function of $N$.
In order to state our main result we first associate, to each $p \in\cU,$ a certain Lyapunov exponent.
Given a symmetric positive definite $m \times m$ real matrix $A,$ and $\beta\in\RR$, define the probability measure $\mu_{A,\beta}$  on $\SS^{m-1}$,  the unit sphere in $\RR^m$, via
\bqn{mAb}
\fo \theta\in\SS^{m-1},\qquad \mu_{A,\beta}(d\theta)&=&\f1{Z(A,1+\beta)}\lan \theta,A\theta\ran^{-1-\beta}\, \sigma(d\theta)\eqn
where, $\sigma$ is the uniform probability measure on $\SS^{m-1},$ $\lan \, \cdot\,,\,\cdot \, \ran$ the Euclidean dot product (not to be confused with the Riemannian metric on $M$) and  $Z(A,1+\beta)$ is the normalization constant.

Define the {\em $\beta$-average eigenvalue} of $A$ as
\bqn{lyap}
\Lambda(A,\beta) = \int_{\SS^{m-1}}\lan \theta,A\theta \ran \mu_{A,\beta}  (d\theta) =
 \frac{Z(A,\beta)}{Z(A,1+\beta)}.
\eqn
Let $\lambda_1(A) \leq \ldots \leq \lambda_m(A)$ be the eigenvalues of $A.$
Observe that $\Lambda(A,\beta)$ only depends on these eigenvalues, because $\sigma$ is invariant by orthogonal transformations and $A$ is orthogonally conjugate to a diagonal matrix. Observe also that
\bqn{estimlyap}
\lambda_1(A) \leq \Lambda(A,\beta) \leq \lambda_m(A).
 \eqn
\begin{rem}\label{rem4}
Inequalities (\ref{estimlyap}) are strict, except when $\lambda_1(A) = \lambda_m(A).$ Furthermore it can be shown (see the appendix Section \ref{5.2})  that for all numbers
$\lambda_- < \lambda < \lambda_+,$ there exists, for $m$ sufficiently large, a $m \times m$ definite positive matrix $A$ such that
$\lambda_1(A) = \lambda_-, \Lambda(A,\beta) = \lambda$ and $\lambda_m(A) = \lambda_+.$
\end{rem}
Given $p \in \cU,$ we let $A_p$ denote the diagonal matrix whose entries $0 < \lambda_1(p) \leq \ldots \leq \lambda_m(p)$ are
the eigenvalues of the Hessian of $U$ at $p.$
Set
 $$\beta_0 \df \frac{m}{2}-1.$$ Our main result is the following.
\begin{theo}\label{theo1}
Let $x \in N$ and $\beta \in \RR.$
\begin{description}
\item[(i)] If $\beta > \beta_0,$ then
 $$
\sum_{p \in \cU} \PP_{x}\lt[ \limsup_{t\ri+\iy} \frac{\ln(d(X\eb(t),p))}{t} \leq - \Lambda(A_p,\beta) (\beta -\beta_0) \rt]\ =\ 1,
 $$
where  each   term in the above sum is positive.
  \item[(ii)] If $\beta < \beta_0,$ then $X\eb$ has, on $N,$ a unique invariant probability distribution  given by
$$\pi_{\beta}(dx) \df \f1{C_{\beta}} U(x)^{-1-\beta}(x) \ell(dx),$$
where $C_{\beta}$ is a normalization constant and $\ell(dx)$ stands for the Riemannian measure. Furthermore:
\begin{description}
\item[(a)] $X\eb$ is positive recurrent on $N,$ meaning that for all $f \in L^1(\pi_{\beta}),$   $\PP_{x}$ a.s.,
$$\lim_{t\ri+\iy} \frac{1}{t} \int_0^t  f(X\eb(s)) ds = \pi_{\beta}(f)$$
\item[(b)]
  There exist positive constants $a,b, \chi$ (depending on $\beta$)  with $\chi < \beta_0 - \beta,$ such that for all $f : N \rightarrow \RR,$ measurable,
$$|\EE_{x} [f(X\eb(t))] -  \pi_{\beta}(f)| \leq \frac{a  e^{-bt}}{d(x,\cU)^{\chi}} \|f\|_{\chi},$$
where
$$\|f\|_{\chi} := \sup_{x \in  N} |f(x)| d(x,\cU)^{\chi}.$$
\end{description}
  \item[(iii)] If $\beta = \beta_0,$ then, for  every neighborhood $O$ of $\cU$, $\PP_{x}$ a.s.,
  $$\lim_{t\ri+\iy} \frac{1}{t} \int_0^t \mathbf{1}_{ \{X\eb(s) \in O\}}ds = 1$$
   \end{description}
 \end{theo}\par
\par
\begin{rem}\label{rem3}
Theorem \ref{theo1} is an improvement over the results of \cite{miclo:hal-04094950}, which showed the a.s.\ convergence of $X\eb$ toward elements of $\cU$ (each being approached with a positive probability) only for $\beta>\beta_\vee \geq \beta_0$ with
\bqn{bv} \beta_\vee &\df&\max_{p \in \cU}\f{\sum_{l\in\lin m\rin}\lambda_l(p)}{2\lambda_{1}(p)}-1,\eqn
and the a.s.\ non-convergence of $X\eb$ toward elements of $\cU$ for $\beta<\beta_\wedge \leq \beta_0$, with
\bqn{bw}
\beta_\wedge &\df& \min_{p \in \cU}\f{\sum_{l\in\lin m\rin}\lambda_l(p)}{2\lambda_{m}(p)}-1.\eqn\par
\end{rem}
\par
\begin{rem}
Here we restrict our attention to dimensions $m\geq 2,$ so that $N $ is connected. The case  $m=1$ which corresponds to the circle is already  treated in \cite{miclo:hal-04094950}.
\end{rem}
\par
\begin{rem}
By Theorem \ref{theo1}, the diffusion $X\eb$ on $N$ is transient for $\beta > \beta_0$ and positive recurrent if and only if $\beta < \beta_0$, due to the fact that $\int_NU^{-1-\beta}\,d\ell=+\iy$ for $\beta\geq \beta_0$.
By standard results (see e.g.\ Kliemann \cite{MR0885138}, Theorem 3.2 applied with $C = N$),
it is then either null recurrent or transient for $\beta = \beta_0.$ It would be interesting to investigate  this situation.
\end{rem}

\section{Euclidean computations}\label{sec:Euclide}
This section considers a situation where the state space $M$ is no longer a compact manifold but the Euclidean space $\RR^m$, with $m\geq 2.$ We state a theorem (Theorem \ref{pro2} below) analogous to Theorem \ref{theo1} (i).  This result is interesting in itself, and its proof allows us to explain, in a simple framework, how to characterize the attractiveness/repulsivity of a global minimum. The main idea is to expand a critical  point to a sphere, using polar decompositions, following \cite{MR3910006}.

Let $U\st \RR^m\ri \RR_+$ be a smooth function with $\min U = 0.$ We assume that for each $p \in \cU \df U^{-1}(0),$ $\Hess U(p)$  is positive definite. In particular, points in $\cU$ are isolated and $\cU$ is therefore countable.

For any fixed $\beta\in\RR$, as in \eqref{Lb}, we are interested in the operator $L_\beta$ defined on $\cC^2(\RR^m),$  via
\bqn{Lbr}
\fo x\in \RR^m,\qquad L_\beta[f](x)&\df& U(x)\tr f(x)-\beta\lan \na U,\na f\ran(x)\eqn
where $\tr$, $\lan\cdot,\cdot\ran$ and $\na,$ respectively denote, the Euclidean Laplacian, scalar product and gradient. Throughout all this section $\|x\| = \sqrt{\lan x, x \ran}$ denotes the Euclidean norm of $x.$
\par
Associated to (\ref{Lbr}) is the  stochastic differential equation
\bqn{sdeRRm}
dX\eb(t) = - \beta \na U(X\eb(t)) dt + \sqrt{2 U(X\eb(t))} dB_t
\eqn
where $B = (B_t)_{t \geq 0}$ is a standard Brownian motion on $\RR^m.$

By local Lipschitz continuity of $\na U$ and $\sqrt{U},$ there exists, for each $x \in \RR^m,$ a unique solution $X\eb : [0, \tau^{\infty}) \rightarrow \RR^m$ starting from $x,$ (i.e.\ $X\eb(0) = x$). Here, $0 < \tau^{\infty} \leq \infty,$ denotes the explosion time of $X\eb$ and
 is characterized by $$\tau^{\infty} > t \Leftrightarrow \|X\eb(t)\| <  \infty.$$ The set $\RR^m \setminus \cU$ is invariant, in the sense that for all $t \geq 0, x \in \RR^m \setminus \cU,$
 $$\PP_x( X\eb(t) \in \RR^m \setminus \cU \, | \tau^{\infty} > t) = 1.$$ The proof of this last point is the same as the proof of Proposition
 \ref{prosemigroup} (i) given in the appendix.

\begin{theo}\label{pro2}
\begin{itemize}
\item[(i)] Suppose $\beta > \beta_0.$ Then, for all
 $x\in \RR^m \setminus \cU$ and $p \in \cU,$
\bqn{B+}
\PP_{x}\lt[ \limsup_{t\ri+\iy} \frac{\ln(\|X\eb(t) -p\|)}{t} \leq - \Lambda(A_p,\beta)(\beta-\beta_0) \rt]&>&0 \eqn where
$A_p, \Lambda(A_p,\beta)$ are defined as in Section \ref{sec:main}.
\item[(ii)]  Suppose $\beta > \beta_0,$ and in addition, that there exist positive constants $\alpha, r$ (possibly depending on $\beta$) such that
 \bqn{alinfini} 2 \beta_0 U(x) - \beta \langle \na U(x), x \rangle &\leq& - \alpha \| x \|^2
 \eqn whenever  $\| x \| \geq r.$
 Then, $\cU$ is finite and for all $x \in \RR^m \setminus \cU,$
\bqn{B_+1} \sum_{p \in \cU} \PP_{x}\lt[ \limsup_{t\ri+\iy} \frac{\ln(\|X\eb(t)-p\|)}{t} \leq - \Lambda(A_p,\beta)(\beta-\beta_0);  \tau^{\infty} = \infty \rt] = 1.  \eqn
\item[(iii)]
Suppose $\beta < \beta_0.$ Then, for all $p \in \cU$ and
 $x\in \RR^m \setminus \{p\}$
\bqn{B_}
\PP_{x}\lt[ \lim_{t \rightarrow \infty} X\eb(t) = p\rt]&=&0. \nonumber \eqn
\end{itemize}
\end{theo}
\begin{rem}
The condition (\ref{alinfini}) is given for its simplicity. However, the conclusion (\ref{B_+1}) holds true under the weaker assumption, implied by (\ref{alinfini}) (see Lemma \ref{explo} below), that $X\eb$ almost surely never explodes (i.e $\tau^{\infty} = \infty$) and eventually enters a ball $B(0,r)$ containing $\cU$ for some $r > 0.$
\end{rem}
The remainder of this section is devoted to the proof of Theorem \ref{pro2}.
We first recall some classical facts about diffusion operators, see e.g.\ Bakry, Gentil and Ledoux \cite{MR3155209}.
The \textbf{carré du champ} $\Gamma_L$ associated to a Markov generator $L$ defined on an algebra $\cA(L)$ is the bilinear functional
defined on $\cA(L)\times \cA(L)$ via
\bq
\fo f,g\in \cA(L),\qquad \Gamma_L[f,g]&\df& L[fg]-fL[g]-gL[f]\eq
(we will denote $\Gamma_L[f]\df \Gamma_L[f,f]$).
\par
The generator $L$ is said to be of \textbf{diffusion}, if $\cA(L)$ is stable by composition with smooth functions and if we have
\bqn{comp1}
L[\varphi(f)]&=&\varphi'(f)L[f]+\f{\varphi''(f)}{2}\Gamma_L[f]\eqn
for any $f\in\cA(L)$ and any function $\varphi$ smooth on the image of $f$.
\par
In this situation we also have, with the same notations,
\bqn{comp2}
\Gamma_L[\varphi(f)]&=&(\varphi'(f))^2\Gamma_L[f]\eqn 
The Markov generator given in \eqref{Lbr} is  of diffusion with $\cA{(L_{\beta})} = \cC^2(\RR^m)$. The corresponding carré du champ is given by
\bqn{G1}
\fo f\in \cC^2(\RR^m),\qquad \Gamma_{L_\beta}[f]&=&2U \|\na f\|^2.
\eqn
Our first goal is to show that, under condition (\ref{alinfini}), $X\eb$ never explodes and always enter the ball $B(0,r).$
For all $s \geq 0,$ we let $$\tau_s = \inf \{t \geq 0 \:  : \|X\eb(t)\| \leq s\}  \mbox{ and } \tau^{s} = \inf \{t \geq 0 \:  : \|X\eb(t)\| \geq s\}.$$
Note that these stopping times depend on $\beta,$ but to shorten notation we omit this dependance in their definition.
\begin{lem} \label{explo}
Under the condition (\ref{alinfini}), $$\PP_x ( \tau^{\infty} = \infty; \tau_{r} < \infty ) = 1$$ for all $x \in \RR^m$ and $r$ is as in \eqref{alinfini}.
\end{lem}
\proof
   Let $V : \RR^m \rightarrow \RR$ be a smooth function coinciding with $\ln(\|x\|^2)$ for $\|x\| \geq r.$
Using the formulaes (\ref{comp1}) and (\ref{G1}) it comes that, for all $\|x\| \geq r,$
$$L_{\beta}(V)(x) = \f{2}{\|x\|^2} \left( 2\beta_0 U(x) -  \beta \lan \na U(x), x \ran \right) \leq - 2 \alpha.$$
In particular, for all $x \in \RR^m, L_{\beta}(V)(x) \leq C$ where $C = \sup_{\{x \in \RR^m \: : \|x\| \leq r\}} |L_{\beta}(V)(x)|.$ Thus, by Ito's formulae, for all $k \geq 1,$
$$\ln(k^2) \PP_x( \tau^k \leq t) \leq \EE_x(V(X\eb(t \wedge \tau^k)) =   V(x) +\EE_x\lt[ \int_0^{t \wedge \tau^k} L_{\beta}[V](X\eb(s)) ds\rt]
\leq V(x) + t C.$$ This shows that $\PP_x( \tau^k \leq t) \rightarrow 0,$ as $k \rightarrow \infty.$ Hence
$\PP_x(\tau^{\infty} < \infty) = 0.$

Now, by Ito formulae again, the process $(M_t)_{t \geq 0}$ defined as
$$M_t := V(X\eb(t \wedge \tau_r))-\ln(r^2) - \int_0^{t \wedge \tau_r} L_{\beta} V(X\eb(s)) ds \geq 2  \alpha (t \wedge \tau_r)$$ is a nonnegative  $\PP_x$ local martingale.
 A nonnegative local martingale may not be a martingale but is always  a supermartingale (Le Gall \cite{MR3497465}, Proposition 4.7). Thus
 $2\alpha \EE_x(t \wedge \tau_r) \leq \EE_x(M_t) \leq V(x)-\ln(r^2) $. Hence $\EE_x(\tau_r) < \infty.$ \wwtbp
 Our next goal is to investigate the behavior of $X\eb$ around a critical point $p \in \cU.$ Without loss of generality, we assume that $p = \{0\}.$ We let $A = \Hess U(0).$ Fix $\epsilon \in(0,1)$ small enough so  that $\cU \cap B(0,\epsilon) = \{0\}.$
 Write any $x\in B(0,\epsilon)\setminus\{0\}$ under its
polar decomposition $x=\rho \theta$ with $\rho\in(0,\epsilon)$ and $\theta \in \SS^{m-1}$.
This decomposition induces the mapping
\bq
P \st \cC^2(B(0,\epsilon))\ni f&\mapsto & P[f]\in \cC^2((0,\epsilon)\times \SS^{m-1})\eq
with
\bqn{P}
\fo (\rho,\theta)\in (0,\epsilon)\times \SS^{m-1},\qquad P[f](\rho,\theta)&\df& f(\rho \theta)\eqn
\par
Endow $\SS^{m-1}$ with its usual Riemannian structure, inherited from $\RR^m$, and denote  $\lan\,\cdot\,,\,\cdot\,\ran_\theta$, $\na_\theta,$ $\mathsf{div}_{\theta}$ and $\tr_\theta$  the corresponding scalar product, gradient, divergence, and Laplace-Beltrami operator. Note that $\lan\,\cdot\,,\,\cdot\,\ran_\theta$ is just the restriction of $\lan\,\cdot\,,\,\cdot\,\ran$ to the tangent space of $\SS^{m-1}$ at $\theta$.
\par
Classical computations in polar coordinates show that for any $f,g\in \cC^2(B(0,\epsilon))$, we have on $(0,\epsilon)\times \SS^{m-1}$,
\bq
P[\lan \na f,\na g\ran]&=&\pa_\rho P[f]\pa_\rho P[g]+\f1{\rho^2}\lan \na_\theta P[f],\na_\theta P[g]\ran_\theta,\\
P[\tr f]&=&\pa^2_\rho P[f]+\f{m-1}{\rho}\pa_\rho P[f]+\f1{\rho^2}\tr_\theta P[f].
\eq
It leads us to introduce the operator $\sL_\beta$ on $\cC^2((0,\epsilon)\times \SS^{m-1})$ defined by
\bqn{sLb}
\sL_\beta\,\cdot\,&\df& \sU\lt(\pa^2_\rho\,\cdot\,+\f{m-1}{\rho}\pa_\rho \,\cdot\,+\f1{\rho^2}\tr_\theta\,\cdot\,\rt)-\beta\lt((\pa_\rho\sU)\pa_\rho\,\cdot\, +\f1{\rho^2}\lan \na_\theta\sU,\na_\theta \,\cdot\,\ran_\theta\rt)
\eqn
where $\sU\df P[U]$. Indeed,  on $\cC^2(B(0,\epsilon))$, we have the intertwining relation
\bq
\sL_\beta\circ P&=&P\circ L_\beta.\eq
\begin{lem} \label{Lextends}
The operator $\sL_\beta$ extends to a diffusion operator, still denoted $\sL_\beta,$ on $\cC^2([0,\epsilon)\times \SS^{m-1})$, whose associated diffusion process $\sX\eb$ leave $\{0\}\times \SS^{m-1}$ invariant.
On $\{0\}\times \SS^{m-1}$, identified with $\SS^{m-1}$,  $\sX\eb$ admits for generator the operator $G_\beta$ acting on $\cC^2(\SS^{m-1})$ via
\bqn{Gb}
\fo \ff\in\cC^2(\SS^{m-1}),\, \qquad G_\beta[\ff] &\df&   \f1{2} \Psi_A^{1+\beta} \mathsf{div}_{\theta} ( \Psi_A^{-\beta} \na_{\theta} f) \eqn where
$\fo \theta\in\SS^{m-1}, \Psi_A(\theta) =  \lan \theta,A\theta\ran.$
Furthermore,  $G_\beta$ has a unique invariant probability measure
on $\SS^{m-1}$, given by $\mu_{A,\beta}$ (see Equation \eqref{mAb}).
\end{lem}
\proof
Our assumptions on $U$ imply that, uniformly over $\theta\in\SS^{m-1}$,
\bqn{lim1}
\lim_{\rho\ri0_+} \f{\sU(\rho,\theta)}{\rho^2}&=&\f12\lan \theta,A\theta\ran,
\eqn
\bqn{lim2}
\lim_{\rho\ri0_+} \f{\pa_\rho\sU(\rho,\theta)}{\rho}&=&\lan \theta,A\theta\ran,
\eqn
\bqn{lim3}
\lim_{\rho\ri0_+} \f{\na_\theta\sU(\rho,\theta)}{\rho^2}&=&A\theta-\lan \theta,A\theta\ran\theta.\eqn
Indeed, by the usual expansion of $U$ around 0, we have
\bq
U(x)&=&U(0)+\lan\na U(0),x\ran+\f12\lan x,\Hess U(0)x\ran+\circ(\lan x,x\ran)
\\&=&
\f12\lan x, Ax\ran +\circ(\lan x,x\ran)\eq
which translates into
\bq
\sU(\rho,\theta)&=&\f{\rho^2}{2}\lan \theta, A\theta\ran+\circ (\rho^2)\eq
leading to the first announced limit (\refeq{lim1}).
Similarly,
\bq
\na U(x)&=& \na U(0)+\Hess U(0) x+\circ (\sqrt{\lan x,x\ran})\\
&=&Ax +\circ (\sqrt{\lan x,x\ran}).\eq
At $x=\rho \theta$ with $\rho>0$, $\pa_\rho\sU(\rho,\theta)\theta$ is the radial part of $\na U(x)$ and
$\na_\theta\sU(\rho,\theta)/\rho$ is the tangential part.
It follows that
\bq
\pa_\rho\sU(\rho,\theta)&=&\lan \na U(x),\theta\ran,\\
\f{\na_\theta\sU(\rho,\theta)}{\rho}&=&\na U(x)-\pa_\rho\sU(\rho,\theta)\theta,\eq
and we get
\bq
\f{\pa_\rho\sU(\rho,\theta)}{\rho}&=&\lan \theta,A\theta\ran+\circ(1),\\
\f{\na_\theta\sU(\rho,\theta)}{\rho^2}&=&A\theta-\lan \theta,A\theta\ran\theta+\circ(1),\eq
leading to the wanted second and third results (\refeq{lim2}) and (\refeq{lim3}).

It follows that for any $F\in\cC^2([0,\epsilon)\times \SS^{m-1})$, we have, uniformly over $\theta\in\SS^{m-1}$,
\bqn{fa1}
\lim_{\rho\ri0_+}\sL_\beta[F](\rho,\theta)&=&\f12\lan \theta,A\theta\ran\tr_\theta F(0,\theta)-\beta \lan A\theta-\lan \theta,A\theta\ran\theta,\na_\theta F(0,\theta)\ran_\theta\eqn
Denoting $\sL_\beta[F](0,\theta)$ the r.h.s.\ enables us to see $\sL_\beta$ as a diffusion operator on $[0,\epsilon)\times \SS^{m-1}$, whose associated diffusion process $\sX\eb$ leaves $\{0\}\times \SS^{m-1}$ invariant, and such that on $\{0\}\times \SS^{m-1}$, identified with $\SS^{m-1}$, its generator  coincides with the operator defined by
$$G_\beta(\ff) := \lan \theta,A\theta\ran\lt(\f12\tr_\theta \ff(\theta)-\beta \lan  b(\theta),\na_\theta \ff\ran_\theta\rt)$$
where
\bq
\fo \theta\in\SS^{m-1},\qquad b(\theta)&\df& \f{A\theta-\lan \theta,A\theta\ran\theta}{\lan \theta,A\theta\ran} =  \f12\na_\theta \ln(\lan \theta,A\theta\ran).\eq It is easily checked that $G_{\beta}$ can be rewritten under the divergence form given by (\refeq{Gb}).
This divergence form implies  that the probability measure $\mu_{A,\beta}$ defined in \eqref{mAb} is invariant. By ellipticity of $G_\beta$ there is no other invariant probability measure.
\par\me \wwtbp
 \begin{lem}
 \label{localconv}
 Suppose $\beta > \beta_0$ and $0 < \lambda < \Lambda(A,\beta).$ There exists $0 < \epsilon_0 \leq \epsilon$ with the property  that for all  $0 < \eta \leq  1,$ there exists $0 < \epsilon_1 < \epsilon_0$ such that for all  $\|x\| \leq \epsilon_1,$
 \bqn{locconv}\PP_{x}\lt[ \limsup_{t\ri+\iy} \frac{\ln(\|X\eb(t)\|)}{t} \leq - \lambda(\beta-\beta_0); \tau^{\epsilon_0} = \infty \rt] = \PP_x \lt[ \tau^{\epsilon_0} = \infty \rt] \geq 1- \eta. \eqn
 If now, $\beta < \beta_0,$ then for all $x \in \RR^m \setminus \{0\},$ $$\PP_{x}\lt[ \lim_{t\ri+\iy} \|X\eb(t)\| = 0 \rt] = 0.$$
 \end{lem}
 \proof The proof follows from the stochastic persistence approach used in \cite{MR3910006}, \cite{Ben18}. Let $\sV$ be the function defined on $(0,\epsilon)\times \SS^{m-1}$
via
\bqn{Vlnr}
\sV(\rho,\theta)&\df&-\ln(\rho).\eqn
We claim that:
\begin{description}
\item{[a]} $\sL_\beta[\sV]$ can be extended into a continuous function $\sH_\beta$ on
$[0,\epsilon)\times \SS^{m-1};$
\item{[b]} $\Gamma_{\sL_\beta}[\sV]$ is bounded on $(0,\epsilon)\times \SS^{m-1};$ and
\item{[c]} $\mu_{A,\beta}[\sH_\beta(0,\cdot)] =  \Lambda(A,\beta) (\beta-\beta_0).$
\end{description}
 Using the form of $\sL_\beta$ (equation (\refeq{sLb})) and the equalities (\refeq{lim1}), (\refeq{lim2}), $(a)$ holds true with \bqn{Hb0}\sH_\beta(0,\theta) &=& (\beta-\beta_0) \lan \theta, A \theta \ran\eqn and $(c)$ directly follows from the definition of $\Lambda(A,\beta).$ For $(b)$, the definition of ${\sL_\beta}$ and $\Gamma_{\sL_\beta},$ lead to \bq
\fo f\in \cC^2((0,\epsilon)\times \SS^{m-1}),\qquad \Gamma_{\sL_\beta}[f]&=&2\sU \lt((\pa_\rho f)^2+\f1{\rho^2}\vert\na_\theta f\vert^2\rt).
\eq
Thus, $$\Gamma_{\sL_\beta}[\sV] = 2 \frac{\sU(\rho,\theta)}{\rho^2}$$ which is bounded in view of (\refeq{lim1}). This concludes the proof of the claim.   \par
If $\beta > \beta_0, \mu_{A,\beta}[\sH_\beta(0,\cdot)] > 0$ and the first assertion of the lemma  follows from Theorem 5.4  in \cite{MR3910006} (to be more precise, this follows from the proof of Theorem~5.4 in \cite{MR3910006}, because  the formulation of Theorem~5.4 in \cite{MR3910006} doesn't specify that $\eta$ can be chosen arbitrary close to one).
If $\beta < \beta_0, \mu_{A,\beta}[\sH_\beta(0,\cdot)] < 0$ and (see e.g.\ \cite{Ben18}, Proposition 8.1 or the proof of Theorem 3.2 (iii) in \cite{MR3910006}) there exist positive constants $\epsilon_1 \leq \epsilon, C$ such that  $\EE_x(\tau^{\epsilon_1}) \leq C \vert\ln(\|x\|)\vert < \infty$ for $x \in B(0,\epsilon_1) \setminus \{0\}.$
   \par\me \wwtbp
 We can now conclude the proof of Theorem \ref{pro2}. We start with assertion $(ii).$  Fix $\beta > \beta_0.$
 For $n \in \NN$ sufficiently large (so that $\Lambda(A_p,\beta) > \f1{n}$) and $p  \in \cU,$ let $\cE_n(p)$ be the event defined as
 $$\cE_n(p) = \left \{\limsup_{t\ri+\iy} \frac{\ln(\|X\eb(t)-p\|)}{t} \leq - (\Lambda(A_p,\beta) -\f1{n}) (\beta-\beta_0)\right \},$$ and let
 $$\cE_n = \bigcup_{p \in \cU} \cE_n(p).$$
The set $\cU$ is finite, since \eqref{alinfini} cannot be satisfied by a point $x\in\cU$ and by consequence $\cU$ is included into the compact ball centered at 0 and of radius $r$. Thus there exists, by Lemma \ref{localconv},  $\epsilon_1 > 0$ such that $$\PP_x(\cE_n(p)) \geq \f1{2}$$ for all  $x \in B(p,\epsilon_1)$ and all  $p \in \cU.$   Let $$\cU_{\epsilon_1} \df \bigcup_{p \in \cU} B(p, \epsilon_1)$$ and $\tau_{\cU_{\epsilon_1}} \df \inf \{t \geq 0 \: : X\eb(t) \in \cU_{\epsilon_1}\}.$
  By ellipticity of $L_{\beta}$ on $\RR^m \setminus \cU$,  $\cU_{\epsilon_1}$ is open and {\em accessible} from all $x \in \RR^m,$ in the sense that $\PP_x( X\eb(t_x) \in\cU_{\epsilon_1}) > 0$ for some $t_x \geq  0.$ Thus,
 by Feller continuity and compactness of $\overline{B(0,r)},$ there exists $\delta > 0$ such that
 $$\PP_x(\tau_{\cU_{\epsilon_1}} < \infty) \geq \delta$$ for all $x \in \overline{B(0,r)}.$   Combined with Lemma \ref{explo}, this proves that $\PP_x( \tau_{\cU_{\epsilon_1}} < \infty) \geq \delta$ for all $x \in \RR^m.$ Thus,  $$\PP_x(\cE_n) \geq \delta/2$$ for all $x \in \RR^m.$  The strong Markov property, implies that $\PP_x( \cE_n) = 1.$ Hence  $$\PP_x (\bigcap_n \cE_n) = 1.$$ This concludes the proof of $(ii).$

We now pass to the proof of $(i)$. Fix $\beta > \beta_0$ and assume without loss of generality that $p = \{0\}.$ Let $\tilde{U} : \RR^m \rightarrow \RR^+$ be a smooth function which coincides with $x \mapsto U(x)$ on a neighborhood of $0,$ with $x \mapsto \|x\|^2$ for $\|x\| \geq 1,$ and such that $U^{-1}(0) = \{0\}.$ Let $\tilde{X}\eb$ be solution to the stochastic differential equation given by (\ref{sdeRRm}) with $\tilde{U}$  instead of $U$ and $(B_t)_{t \geq 0}$ the  Brownian motion governing $X\eb.$ The process $\tilde{X}\eb$ satisfies (\ref{B_+1}) because $\tilde{U}$ satisfies (\ref{alinfini}). Thus, by Theorem \ref{pro2} $(ii)$ (applied to $\tilde{X}\eb$) and  Lemma  \ref{localconv} (applied to $X\eb$), there exist $0 < \epsilon_1 < \epsilon_0 \leq  \epsilon$ such that whenever $\|x\| \leq \epsilon_1,$
 $$\frac{1}{2} \leq \PP_x \lt[ \tau^{\epsilon_0} = \infty \rt] = \PP_{x,x} \lt[ \limsup_{t\ri+\iy} \frac{\ln(\|\tilde{X}\eb(t)\|)}{t} \leq - \Lambda(A,\beta)(\beta-\beta_0);  \tau^{\epsilon_0}= \infty  \rt]$$
 $$ = \PP_{x} \lt[ \limsup_{t\ri+\iy} \frac{\ln(\|X\eb(t)\|)}{t} \leq - \Lambda(A,\beta)(\beta-\beta_0);  \tau^{\epsilon_0} = \infty \rt].$$ Here $\PP_{x,x}$ stands for the law of $(X\eb,\tilde{X}\eb)$ starting from $(x,x).$ Since $\PP_x(\tau_{\epsilon_1} < \infty) > 0$ for all $x \in \RR^m \setminus \{0\},$ the proof of (\ref{B+}) follows.
\par
Finally (iii) is an immediate consequence of the second part of Lemma \ref{localconv} (recall that 0 was an arbitrary point of $\cU$, up to a translation).

\section{Proof of Theorem \ref{theo1}}\label{poT}
\label{sec:proofs}
\subsection{Proof of Theorem \ref{theo1} $(i)$}
\label{poT1}
The proof is similar to that of Theorem  \ref{pro2}. We begin by proving a Riemannian version of Lemma \ref{localconv}. The proof of Theorem \ref{theo1} $(i)$ will then follow by an argument similar to that given at the end of  Section~\ref{sec:Euclide}.

Let $y \in \cU$ and let $B_M(y,\epsilon)$  be the Riemannian ball with center $y$ and radius $\epsilon,$ where  $\epsilon>0$ is sufficiently small so that
 \begin{itemize}
 \item the only critical point for $U$ in $B_M(y,\epsilon)$ is $y$,
 \item the exponential mapping $\exp_y\st T_yM\ri M$ is a diffeomorphism between the tangent ball $B(0,\epsilon)$ of $T_yM$ and $B_M(y,\epsilon)$.
 \end{itemize}
 \par
Recall that the exponential mapping $\exp_{y}\st T_{y}M\ri M$  associates to any tangent vector $v\in T_{y}M$ the point $x\in M$ which is the position at time 1 of the (constant speed) geodesic starting at time 0 from $y$ with speed $v$.

Consider $(e_1, e_2, ..., e_m)$ an orthonormal basis of $T_yM$ consisting of eigenvectors associated to the eigenvalues $(\lambda_1, \lambda_2, ..., \lambda_m)$ of the Hessian of $U$ at the critical point $y$. A priori this Hessian is a bilinear form on $T_yM$, but the Euclidean structure of $T_yM$ enables us to see it as a symmetric endomorphism on $T_yM$, and  $(\lambda_1, \lambda_2, ..., \lambda_m)$ and $(e_1, e_2, ..., e_m)$ correspond to its spectral decomposition.\par
Let $(v_1,v_2, ..., v_m)$ be the coordinate system associated to $(e_1, e_2, ...,e_m)$ on $B(0,\epsilon)$.
Such a coordinate system based on the exponential mapping is said to be a normal. From now on and until the end of this section, we identify a map $f : B_M(y,\epsilon) \rightarrow \RR$ with $f \circ \exp_y : B(0,\epsilon)\rightarrow \RR,$ and write $f(v)$  for $f \circ \exp_y(v).$    Under this identification, the matrix corresponding to the Hessian at $y$ admits the classical form
\bq
(\pa_{k,l} U(0))_{k,l\in \lin m\rin}\eq
where $\pa_k$ is a shorthand for $\f{\pa}{\pa_{v_k}}$
The introduction of the lecture notes of Pennec \cite{Pennec_geometric} is a convenient reference for these assertions (a more thorough exposition can be found in the book of Gallot,  Hulin and Lafontaine \cite{MR2088027}).\par
A first interest of the normal coordinate system $(v_1,v_2, ..., v_m)$ on $B(0,\epsilon)$ is that we can consider the corresponding polar decomposition as in the previous section: each $v=(v_1, v_2, ..., v_m)\in B(0,\epsilon)\setminus\{0\}$ can be uniquely written under the form $\rho \theta$ with $\rho\in(0,\epsilon)$ and $\theta\in\SS^{m-1}$, where the basis $(e_1, e_2, ..., e_m)$ enables us to identify $T_yM$ with $\RR^m$.
\par
Before going further, let us recall some other traditional notations and facts from Riemannian geometry.
For any $v\in B(0,\epsilon)$,  denote  $g(v)\df(g_{k,l}(v))_{k,l\in \lin m\rin}$ the matrix of the pull-back of the Riemannian metric:
for any vectors $b$ and $\wi b$ from $T_{\exp_y(v)}M$, identified with their coordinates $(b_k)_{k\in\lin m\rin}$ and $(\wi b_k)_{k\in\lin m\rin}$ in the basis $(\pa_k)_{k\in\lin m\rin}$,
we have
\bq
\lan b,\wi b\ran_{v}&=&\sum_{k,l\in\lin m\rin} g_{k,l}(v)b_k\wi b_l\eq
\par
The determinant of $g(v)$ and the inverse matrix $g^{-1}(v)$ are respectively denoted $\vert g\vert(v)$ and \linebreak $(g^{k,l}(v))_{k,l\in \lin m\rin}$.
For any  smooth function $f,$ the expressions of its gradient and  Laplacian are given by
\bq
\na f(v)&=&\lt(\sum_{l\in\lin m\rin}g^{k,l}(v)\pa_l f(v)\rt)_{k\in\lin m\rin}\\
\tr f(v)&=&\f1{\sqrt{\vert g\vert (v)}}\sum_{k,l\in\lin m\rin} \pa_k\lt(\sqrt{\vert g\vert }g^{k,l}\pa_l f\rt)(v)\\
&=&\sum_{k,l\in\lin m\rin}g^{k,l}(v)\lt(\pa_{k,l} f(v)-\sum_{j\in\lin m\rin}\Gamma^{j}_{k,l}(v)\pa_j f(v)\rt)
\eq
where $\Gamma^{j}_{k,l}(v)$ are the Christoffel symbols at $v$, see for instance the listing \cite{enwiki:1157164633} (again we abuse notation in the r.h.s\ by identifying $f$ with its formulation in the coordinate system $v=(v_1,v_2, ..., v_m)$).
\par
A second interest of the normal coordinate system is that at 0, we recover the usual notions: $g(0)$ is the identity matrix and the Christoffel symbols all vanishes at $0$.\par
The above expressions lead to the following formulation of the generator $L_\beta$ defined in \eqref{Lb}:
\bqn{Lb2}
L_\beta\,\cdot\,&=& U\sum_{k,l\in\lin m\rin}g^{k,l}\lt(\pa_{k,l}\,\cdot\,-\sum_{j\in\lin m\rin}\Gamma^{j}_{k,l}\pa_j \,\cdot\,\rt)
-\beta \sum_{k,l\in\lin m\rin}g^{k,l}\pa_kU\pa_l \,\cdot\,
\eqn
Again we are slightly abusing notations by calling it $L_\beta$ too, especially as we see it as  only defined  on $\cC^2(B(0,\epsilon))$.
\par
The associate carré du champ is given as
\bqn{GLb2}
\Gamma_{L_{\beta}}\,\cdot\,&=&  2 U  \sum_{k,l\in\lin m\rin}g^{k,l} \pa_k \,\cdot\, \pa_l \,\cdot\,
\eqn
(this is a consequence of the algebraic relation $\Gamma_{\pa_k\pa_l}\,\cdot\,=2\pa_k \,\cdot\, \pa_l \,\cdot\,$, even if $\pa_k\pa_l$ is not a Markov generator, i.e.\ when $k\neq l$).\par
Consider  the mapping $P$ associated  in \eqref{P} to the polar decomposition. Since $P$ is invertible from $\cC^2(B(0,\epsilon))$
to $\cC^2((0,\epsilon)\times\SS^{m-1})$, there is a unique diffusion generator $\sL_\beta$ acting on $\cC^2((0,\epsilon)\times\SS^{m-1})$ such that
\bq
\sL_\beta\circ P&=&P\circ L_\beta\eq
\par
To compute $\sL_\beta$, let us write that for any $v\in B(0,\epsilon)\setminus\{0\}$,
\bq
\rho&=&\sqrt{\sum_{k\in\lin m\rin}v_k^2}\\
\fo l\in \lin m\rin,\qquad \theta_l&=&\f{v_l}{\rho}\eq
\par
It follows that for any $k\in \lin m\rin$,
\bq
\pa_k\rho&=&\f{v_k}{\rho}\ =\ \theta_k\\
\fo l\in\lin m\rin,\qquad \pa_k\theta_l&=&\f{\delta_{k,l}}{\rho}-\f{v_l}{\rho^2}\pa_k\rho\ =\ \f1\rho(\delta_{k,l}-\theta_k\theta_l)\eq
where $\delta_{k,l}$ is the Kronecker symbol.\par
It follows that
\bqn{ppp}
\pa_k&=&\theta_k\pa_\rho+\f1\rho\sum_{l\in\lin m\rin}(\delta_{k,l}-\theta_k\theta_l)\pa_{\theta_l}\eqn
and by composition, for any $k,l\in\lin m\rin$, we can also write $\pa_{k,l}$ in terms of $\pa_\rho$, $\pa_\rho^2$, $\pa_{\theta_{i}}$ and $\pa_{\theta_{i},\theta_{j}}$,
for $i,j\in \lin m\rin$. Replacing these expressions in \eqref{Lb2}, we get the formula for $\sL_\beta$ in terms of differentiations of order 1 and 2, with respect to $\rho$ and the $\theta_l$, $l\in\lin m\rin$.
\par
In order to apply the general method of  \cite{MR3910006} as in Section \ref{sec:Euclide}, we need to check the three facts respectively listed in the following lemmas.
\begin{lem}\label{l12}
 For any $F\in\cC^2([0,\epsilon)\times \SS^{m-1})$, we have, uniformly over $\theta\in\SS^{m-1}$,
\bq
\lim_{\rho\ri0_+}\sL_\beta[F](\rho,\theta)&=&G_\beta[F(0,\cdot)](\theta)\eq
where $G_\beta$ is given in \eqref{Gb}.
\end{lem}
\proof
For any $v\in B(0,\epsilon)$, define
\bq
\fo k,l\in \lin m\rin,\qquad \wi g^{k,l}(v)&\df& \delta_{k,l}\\
\fo j, k,l\in \lin m\rin,\qquad \wi\Gamma_{k,l}^j(v)&\df& 0\eq
and in analogy with \eqref{Lb2},
\bq
\wi L_\beta\,\cdot\,&=& U\sum_{k,l\in\lin m\rin}\wi g^{k,l}\lt(\pa_{k,l}\,\cdot\,-\sum_{j\in\lin m\rin}\wi\Gamma^{j}_{k,l}\pa_j \,\cdot\,\rt)
-\beta \sum_{k,l\in\lin m\rin}\wi g^{k,l}\pa_kU\pa_l \,\cdot\,
\eq
\par
This operator coincides with the restriction of \eqref{Lbr} to $B(0,\epsilon)$. It follows from \eqref{fa1} that uniformly over $\theta\in\SS^{m-1}$,
\bq
\lim_{\rho\ri0_+}\wi\sL_\beta[F](\rho,\theta)&=&G_\beta[F(0,\cdot)](\theta)\eq
where the operator $\wi\sL_\beta$ is such that $\wi\sL_\beta\circ P=P\circ\wi L_\beta$.
\par
Thus to get the wanted result, it is sufficient to show that
\bqn{LL}
\lim_{\rho\ri0_+}(\sL_\beta-\wi\sL_\beta)[F](\rho,\theta)&=&0\eqn
\par
This convergence is a consequence of the writing
\bq
(\sL_\beta-\wi \sL_\beta)[F]&=& U\sum_{k,l\in\lin m\rin}(g^{k,l}-\wi g^{k,l})\lt(\pa_{k,l}F-\sum_{j\in\lin m\rin}\Gamma^{j}_{k,l}\pa_j F\rt)-
U\sum_{k,l,j\in\lin m\rin}\wi g^{k,l}(\Gamma^{j}_{k,l}-\wi\Gamma^{j}_{k,l})\pa_j F
\\&&-\beta \sum_{k,l\in\lin m\rin}(g^{k,l}-\wi g^{k,l})\pa_kU\pa_l F
\eq
and of the following facts, valid uniformly in $\theta\in\SS^{m-1}$ as $\rho$ goes to $0_+$:
\begin{itemize}
\item According to \eqref{ppp}, for any $k,l\in\lin m\rin$, $\pa_k F$ is of order $1/\rho$ and $\pa_{k,l} F$ is of order $1/\rho^2$.
\item Due to the regularity of $g$ and of the Christoffel symbols, for any $k,l\in\lin m\rin$, $g^{k,l}-\wi g^{k,l}$ and $\Gamma^j_{k,l}-\wi \Gamma^j_{k,l}$ are of order $\rho$.
\item By the assumption that $y$ is a global minimal, $U$ is of order $\rho^2$ and $\pa_k U$ is of order $\rho$, for any $k\in\lin m\rin$.
\end{itemize}
\wwtbp
\par
We have seen in the previous section that $G_\beta$ is reversible with respect to the probability measure $\mu_{A,\beta}$ defined  in \eqref{mAb}, where here $A \df A_y$ is the diagonal matrix whose entries are the eigenvalues of the Hessian of $U$ at $y\in \cU$.
To  continue the method of  \cite{MR3910006}, we also need the two following ingredients.
\begin{lem}
\label{GLV} Consider the function
$\sV$ defined on $(0,\epsilon)\times \SS^{m-1}$
via
\bq
\sV(\rho,\theta)\ \df\ -\ln(\rho).\eq
The function
 $\Gamma_{\sL_\beta}[\sV]$ is bounded on $(0,\epsilon)\times \SS^{m-1}$ and the function $\sL_\beta[\sV]$ can be extended into a continuous function $\sH_{\beta}$ on
$[0,\epsilon)\times \SS^{m-1}$ satisfying \eqref{Hb0} and thus $\mu_{A,\beta}[\sH_\beta(0,\cdot)] =  \Lambda(A,\beta) (\beta-\beta_0).$
 \end{lem}
 \proof
We have $\Gamma_{\sL_\beta}[\sV]=P[\Gamma_{L_\beta}[V]]$ with $V(v)=-\f1{2}\ln(\sum_{k\in\lin m \rin} v_k^2)$, so it is sufficient to see that
$ \Gamma_{L_\beta}[V]$ is bounded on $B_M(y,\epsilon)\setminus\{y\}$.
Expanding $U(v)$ near $0$ in the normal coordinate system $v=(v_1, v_2, ..., v_d)$, we get for $v$ small
\bq
U(v)&\sim& \f12\sum_{l\in\lin m\rin} \lambda_l v_l^2
\eq
Hence, using (\ref{GLb2}),  $$\Gamma_{L_\beta}[V](v) \sim  \frac{\sum_{l\in\lin m\rin} \lambda_l v_l^2}{ (\sum_{l\in\lin m\rin} v_l^2)^2}
 \sum_{k,l\in\lin m\rin} g^{k,l}(0) v_k v_l = \frac{\sum_{l\in\lin m\rin} \lambda_l v_l^2}{\sum_{l\in\lin m\rin} v_l^2}.$$
(see also \cite{miclo:hal-04094950}).
\par
This proves the wanted boundedness.
\par
For the wanted convergence, in view of the computations of the previous section, it is sufficient to see that \eqref{LL} holds with $F$ replaced by $\sV$.
Note that when applied to a function only depending on $\rho$, as $\sV$, \eqref{ppp} reduce to $\pa_k=\theta_k\pa_\rho$.
It follows that $\pa_k\sV$ is of order $1/\rho$ and $\pa_{k,l}\sV$ is of order $1/\rho^2$. This observation enables us to use the same arguments as in the end of the proof of Lemma \ref{l12} to conclude that \eqref{LL} holds with $F$ replaced by $\sV$.
 \wwtbp
 A Riemannian version of Lemma \ref{localconv} follows directly from the preceding lemma, the proof being exactly the same as the proof of Lemma \ref{localconv}. The proof of Theorem \ref{theo1} $(i)$ then follows (almost) verbatim along the lines of the arguments given in the preceding section just after the proof of Lemma \ref{localconv}.

\subsection{Proof of Theorem \ref{theo1} $(ii)$}
 Let $V: N \rightarrow \RR, x \mapsto \ln(U(x)^{-\beta}).$
Observe that for all $f \in C^2(N),$
$$\dive{(e^V \nabla f)} = e^V (\langle \nabla V, \nabla f \rangle + \Delta f) = U^{-\beta -1} L_{\beta} f.$$
Let $C^2_c(N)$ be the set of  $f  \in C^2(N)$  having compact support. Then,
 for all $f \in C^2_c(N),$
$$\int_N L_{\beta} f \, d\ell_{\beta} = 0,$$ where $\ell_{\beta}$ is the measure on $N$ defined as
$$\ell_{\beta}(dx) \df U(x)^{-(1+ \beta)}\ell(dx).$$
Let $p \in \cU.$ By Morse's lemma, there is a smooth chart at $p$ such that, in this chart system, $U$ writes $x \mapsto \|x\|^2 = \sum_{i = 1}^m x_i^2.$ Since the map $x \mapsto \|x\|^{-2(\beta+1)}$ is locally integrable (i.e.\ in a neighborhood of $0_{\RR^m}$) if and only if  $2(\beta + 1) < m,$ it comes that $\int_N U(x)^{-(1+ \beta)} \ell(dx) <  \infty$ if and only if  $2(\beta + 1) < m,$ that is $\beta < \beta_0.$

Assuming   $\beta < \beta_0,$ the probability measure $$\pi_{\beta}(dx) \df \f1{C_{\beta}} \ell_{\beta}(dx)$$  (where $C_{\beta}$ is a normalization constant) satisfies
\bqn{piinv}
\int_N L_{\beta} f \, d\pi_{\beta} = 0
\eqn  for all $f \in C^2_c(N).$
 Observe that there is no evidence that  the set $C^2_c(N)$ is a core for $\cL_{\beta},$ so that we cannot immediately deduce from (\ref{piinv})  that $\pi_{\beta}$ is an invariant probability measure of $X\eb.$ However, by Theorem 9.17 page 248 in Ethier and Kurtz \cite{EK} (originally due to Echeverria \cite{Echeverria}) the following properties $(a)-(d)$ below  ensure that $\pi_{\beta}$ is invariant:
 \begin{itemize}
 \item[(a)]
 The space $N$ is a separable locally compact metric space (for which the space $\hat{C}(N)$ of continuous function "vanishing at infinity" coincide with $\{f \in C^0(M) \: : f|_{\cU} = 0\}$);
 \item[(b)] The set $C^2_c(N)$ is an algebra dense in $\hat{C}(N);$
  \item[(c)] The operator $L_{\beta}:C_c^2(N) \rightarrow \hat{C}(N),$ satisfies the positive maximum principle;
   \item[(d)] The martingale problem for $(L_{\beta}, C_c^2(N))$ is well-posed: for all $x \in N,$  $\PP_x^{\beta}$ (the law of $X\eb$ starting from $X\eb(0) = x$) is the unique probability on $D([0,\infty), N)$ such that
    $f(X(t)) - \int_0^t L_{\beta} f(X(s)) ds$ is a $\PP_x^{\beta}$-martingale and $\PP_x^{\beta}[X(0) =x ] = 1$, where $(X(t))_{t\geq 0}$ is the canonical process on $D([0,\infty), N)$.
     \end{itemize}
     Properties $(a)-(c)$ are easy to verify. Property $(d)$ follows from, 
     on one hand, that for any $\epsilon>0$ sufficiently small, the stopped martingale problem on $N_\epsilon\df\{x\in M\st U(x)\geq \epsilon\}$ is well-posed by uniform ellipticity of  $L\eb$
     on $N_\epsilon$, and on the other hand, that these localized martingale problems can next be extended to the whole state space $N$. 
 For instance, corresponding precise statements are  found in Ethier and Kurtz \cite{EK}, see Theorem 5.4 page 199, providing the existence of a solution of the stopped martingale problem on the $N_\epsilon$, but also of the martingale problem on $N$, Theorem 4.1 page 182 for the uniqueness of stopped martingale problems on the $N_\epsilon$, and Theorem 6.2 page 217, for the deduction of the uniqueness of the solution of the martingale problem on $N$ by localization.
    \par

\paragraph{$\bullet (ii) (a):$} follows from the fact that
a strong Feller process on a connected space having an invariant probability measure with full support, is positive recurrent (see e.g.\ \cite{BH22}, Corollary 7.10 for a statement on discrete time Markov chains and Proposition 4.58 (ii) for the application in continuous time). In particular, it is uniquely ergodic (i.e.\  its invariant probability measure is unique). Here the strong Feller property of $X\eb$  on $N$ follows from Proposition \ref{prosemigroup}.

\paragraph{$\bullet (ii) (b):$} 
The following lemma is a consequence of Lemma \ref{GLV} and the stochastic persistence approach exposed in \cite{Ben18}, \cite{MR3910006}.
\begin{lem}
Assume $\beta < \beta_0.$ Then,
there exist a  continuous map $W: N \rightarrow \RR^+,$ $0 \leq \rho < 1,\, \chi > 0,\, \kappa \geq 0$ and $T > 0$  such that
 \begin{description}
 \item[(i)] $W(x) = d(x,\cU)^{-\chi}$ on a neighborhood of $\cU,$
 \item[(ii)] $P\eb_T W \leq \rho W + \kappa.$
 \end{description}
 \end{lem}
 \proof
  For $y \in \cU,$ and $\epsilon > 0$ sufficiently small, let  $V_y : M \setminus \{y\} \rightarrow \RR^+$ be a smooth map such that 
   $$P[V_y\circ \exp_y](\rho,\theta) = \sV(\rho,\theta):= -\ln(\rho)$$ whenever $\rho < \epsilon,$ where, using the notation of Section \ref{poT1}, $\sV : (0,\epsilon) \times \SS^{m-1} \rightarrow \RR$ is  as in Lemma \ref{GLV} and $P$ is the mapping induced by the polar decomposition as in (\ref{P}). Because  $\Gamma_{\sL_\beta}[\sV]$ is bounded on $(0,\epsilon)\times \SS^{m-1}$ and $\mu_{A,\beta}[\sH_\beta(0,\cdot)] =  \Lambda(A,\beta) (\beta-\beta_0) < 0,$ it is possible, for $\epsilon$ sufficiently small, to find numbers $\chi, T > 0, \kappa$ and $0 \leq \rho < 1$ such that
 $$P\eb_T(e^{\chi V_y}) \leq \rho e^{\chi V_y} + \kappa$$ on  $M \setminus \{y\}$ (see \cite{Ben18}, Proposition 8.2).
The mapping $W : N \rightarrow \RR^+,$ defined as $W(x) = \sum_{y \in \cU} e^{\chi V_y}$  satisfies the conditions of the lemma. \wwtbp
By ellipticity of $L\eb$ on $N$, every point $p \in N$ is an {\em accessible Doeblin} point for $P\eb_T.$ Combined with the preceding lemma this proves assertion $(ii) (b)$ of Theorem \ref{theo1} (see e.g.\ Theorem 8.15 in \cite{BH22}).
\subsection{Proof of Theorem \ref{theo1} $(iii)$}
It follows from compactness of $M$ and Feller continuity of $X\eb$ that, with $\PP_{x}$ probability one, every limit point (for the weak* topology) of the family
$$\left \{\frac{1}{t} \int_0^t \delta_{X\eb_s} ds \right \}_{t \geq 0}$$  is an invariant probability of $X\eb$ (see e.g.\ \cite{BH22}, Theorem 4.20 combined with Propositions 4.57 and 4.58). It then suffices to show that for $\beta = \beta_0,$ every invariant probability of $X^{(\beta_0)}$ is supported by $\cU,$ or equivalently, that every ergodic probability measure of  $X^{(\beta_0)}$ is a Dirac measure $\delta_p$ for some $p \in \cU.$ We proceed by contradiction. Suppose that there exists  an ergodic probability measure of  $X^{(\beta_0)}, \mu$ with $\mu(N) > 0.$ Then $\mu(N) = 1$ (by invariance of $N$) and, by ellipticity of $X^{(\beta_0)}$ on $N,$ $\mu$ is absolutely continuous with respect to $\ell(dx),$ hence  also with respect to $\ell_{\beta_0}(dx).$ That is $\mu(dx) = f(x) \ell_{\beta_0}(dx)$ with $f \geq 0$ measurable and  $\ell_{\beta_0}[f]= 1.$  We claim that $f$ is almost surely constant. This is in contradiction with the fact that $\ell_{\beta_0}(N) =  \infty.$ It  remains to prove the claim. First assume that $\|f\|_{\infty} = \sup_{x \in N } |f(x)| < \infty.$ Then, $f \in L^2(\ell_{\beta_0})$ because $\ell_{\beta_0}[f^2] = \mu[f] \leq \|f\|_{\infty}.$
 Thus, $$\ell_{\beta_0}[ (P^{\beta_0}_t f - f)^2] =  \ell_{\beta_0} [(P^{\beta_0}_t f)^2 + g]$$ where
 $g \df f^2 -2f P^{\beta_0}_t f \in L^1(\ell_{\beta_0})$ and $\ell_{\beta_0}[g] = - \mu[f].$ Thus,
 $$\ell_{\beta_0}[(P^{\beta_0}_t f - f)^2] = \ell_{\beta_0}[(P^{\beta_0}_t f)^2] - \mu(f) = \ell_{\beta_0}[(P^{\beta_0}_t f)^2 - f^2] \leq 0$$ where the last inequality follows from Jensen's inequality. This shows that $\ell_{\beta_0}$-almost surely,  $P^{\beta_0}_t f = f$, and also $\mu$-almost surely. By ergodicity $f$ is $\mu$-almost surely constant. Suppose now that $\|f\|_{\infty} = \infty.$ Set $f_n = \min\{f,n\}$ and $\mu_n(dx) = f_n(x) \ell_{\beta_0}(dx).$ For every Borel set $A \subset N,$
 $$(\mu_n P^{\beta_0}_t)(A) = \mu_n P^{\beta_0}_t (A \cap \{f \leq n\}) + (\mu_n P^{\beta_0}_t)(A \cap \{f > n\}) \leq (\mu P^{\beta_0}_t) (A \cap \{f \leq n\}) + n (\ell_{\beta_0} P^{\beta_0}_t) (A \cap \{f > n\})$$
$$ = \mu (A \cap \{f \leq n\}) + n \ell_{\beta_0}(A \cap \{f > n\}) = \mu_n(A).$$
This shows that $\mu_n$ is excessive, hence invariant because every finite excessive measure is invariant (see e.g.\ \cite{BH22}, Lemma 4.25). By what precedes, $f_n$ is $\mu$-almost surely constant. Thus $f$ is $\mu$-almost surely constant. This concludes the proof of the claim.

\section{Appendix}
\subsection{The diffusion process generated by $L_{\beta}$ and Proposition \ref{prosemigroup}}\label{5.1}
Here we  briefly explain how the diffusion $X\eb$ can be constructed and give a  proof of Proposition  \ref{prosemigroup}.

By Nash's embedding theorem, we can assume without loss of generality that $M$ is a Riemannian submanifold of $\RR^n$ (equipped with its Euclidean scalar product $\langle \, , \, \rangle$) for some $n$ sufficiently large. For reasons that will become clear shortly, we write  $\na_M, \tr_M, \mathsf{div}_M$ the gradient, Laplacian, and divergence on $M,$ and  $\na, \mathsf{div},$  the gradient and divergence on $\RR^n.$
If $F$ is a smooth vector field on $M$ and $\tilde{F}$ a smooth globally integrable vector field on $\RR^n$ such that $\tilde{F}|_M = F,$ then $\tilde{F}$ and  $F$, induce  operators on $C^1(\RR^n)$ and $C^1(M)$ respectively defined by:
 $$\tilde{F}(\tilde{f})(x) = \langle \nabla \tilde{f}(x), \tilde{F}(x) \rangle = \frac{d (\tilde{f} \circ \Psi_t(x))}{dt}|_{t = 0}$$ for all $\tilde{f} \in C^1(\RR^n),$ and  $x \in \RR^n;$
 $$F(f)(x) = \langle \nabla_M f(x), F(x) \rangle = \frac{d (f \circ \Psi_t(x))}{dt}|_{t = 0}$$ for all $f \in C^1(M),$ and $x \in M.$ In both formulae,    $(\Psi^i_t)_{t \in \RR}$ denotes the flow on $\RR^n$ induced by $\tilde{F}.$

 A direct consequence of the right hand side equalities is that
\bqn{opF}
\tilde{F}(\tilde{f})|_M = F(f)
\eqn
for every $f \in C^1(M)$ and $\tilde{f} \in C^1(\RR^n)$ such that $f = \tilde{f}|_M.$

Let $(e_1, \ldots, e_n)$ be the canonical basis of $\RR^n.$ For $i = 1, \ldots, n$ and $x \in M,$ let $E_i(x) \in T_xM$ be the orthogonal projection of $e_i$ onto $T_xM.$ Let $\tilde{E}_i$ be a smooth  vector field on $\RR^n,$ having compact support, such that $\tilde{E}_i|_M = E_i.$ It is not hard to show that such a vector field exists. One can, for example, proceed as follows. Let  $\cM \subset \RR^n$ be a  normal tubular neighborhood of $M.$ Every point $y \in \cM$ writes uniquely $y = x + v$ with $x \in M$ and $v \in T_xM^{\perp}.$ The map $r: \cM \ni  x + v \mapsto x \in M,$ is a smooth retraction. It suffices to set  $\tilde{E}_i(x) = \eta(x) E_i(r(x))$ if $x \in \cM$ and $\tilde{E}_i(x) = 0$ otherwise, where  $0 \leq \eta \leq 1$ is a smooth function with compact support in $\cM$ such that $\eta|_M = 1.$

The following, key property, is proved in Stroock \cite{Stroock}, Section 4.2.1. For the reader's convenience we provide an alternative short proof.
\begin{lem}
\label{DeltaM} For every $f \in C^2(M)$ and $\tilde{f} \in \cC^2(\RR^n),$ such that $f = \tilde{f}|_M,$ one has
$$\sum_{i = 1}^N \tilde{E_i}^2(\tilde{f})|_M = \tr_M(f)$$
\end{lem}
\proof  Let $F$ be a $\cC^1$ vector field on $M,$  and $\tilde{F}$ a $\cC^1$ vector field on $\RR^n$ such that $\tilde{F}|_M = F.$ For all $x \in \RR^n$ $\mathsf{div} \tilde{F}(x)$ equals the trace of the Jacobian matrix  $D\tilde{F}(x),$ while for all $x \in M,$ $\mathsf{div}_M F(x)$ equals the trace of the $d \times d$ matrix $(\langle D\tilde{F}(x) u_i, u_j \rangle)_{i,j}$ where $u_1, \ldots, u_d$ is an (arbitrary) orthonormal basis of $T_xM.$ This has the interesting consequence that
$$\mathsf{div}_M(F) = \mathsf{div}(F \circ r)|_M$$ where $r : {\cM} \rightarrow M$ is the retraction defined above. Let $f \in \cC^2(M).$ Then,
$$\na_M f = \sum_{i = 1}^n \langle \na_M f, e_i \rangle e_i = \sum_{i = 1}^n \langle \na_M f, E_i \rangle e_i = \sum_{i = 1}^n E_i(f)  e_i.$$ Thus,
$$\tr_M f := \mathsf{div}_M (\na_ M f) = \mathsf{div} (\na_M (f) \circ r)|_M = \sum_{i = 1}^n \mathsf{div} [(E_i(f) \circ r)  e_i]|_M $$
$$ = \sum_{i = 1}^n \langle \na (E_i(f) \circ r)|_M, e_i \rangle = \sum_{i = 1}^n \langle \na_M E_i(f), e_i \rangle = \sum_{i = 1}^n E_i^2(f).$$ Here we have used the fact that $\na (f \circ r)|_M = \na_M f$ for all $f \in \cC^1(M).$ \wwtbp
Now, let $\tilde{U} : \RR^n \rightarrow \RR_+$ be a smooth function such that   $\tilde{U}|_M = U,$ $\sqrt{\tilde{U}}$ is Lipschitz and $\na {\tilde{U}}$ has compact support. For instance $\tilde{U}(x) = \eta(x) U(r(x)) + 1-\eta(x)$ for $x \in \cM$ and $\tilde{U}(x) = 1$ otherwise, where $\eta, r$ are as above.
 Here, the Lipschitz continuity of  $\sqrt{\tilde{U}}$ follows from the fact that $r$ is smooth  and that, by assumption, the zeroes of $U$ are non-degenerate.

Consider the stochastic differential equation  on $\RR^n$  defined by
\begin{eqnarray}
\label{sde}
  dX(t) &=&  (-\beta - \f1{2}) \na \tilde{U}(X(t)) dt  \nonumber \\
   &  & +  \sum_{i = 1}^n \left ( \f1{2} \langle \na \tilde{U}(X(t)),\tilde{E}_i(X(t) \rangle  \tilde{E}_i(X(t))  +  \tilde{U}(X(t))  D\tilde{E}_i(X(t))\cdot \tilde{E}_i(X(t) \right ) dt   \nonumber \\
    & & + \, \sqrt{2 \tilde{U}(X(t))} \sum_{i = 1}^n \tilde{E}_i(X(t))  dB^{i}(t)
\end{eqnarray}
where $B = (B^1(t), \ldots, B^n(t))_{t \geq 0}$ is a $n$-dimensional Brownian motion with $B(0) = 0.$

 Since the coefficients of (\ref{sde}) are globally Lipschitz and bounded, the following properties $(a),(b),(c)$ are classical (see e.g.\ Le Gall \cite{MR3497465}, Theorems 8.3 and 8.7 for $(a)$ and $(b),$ and Kunita \cite{MR1070361}, Theorem 4.5.1 for $(c)$) :
\begin{description}
\item[(a)] For all  $x \in \RR^n,$  there is a unique strong solution $\RR_+\ni t \mapsto X\ebx(t)$ to (\ref{sde}) such that $X\ebx(0) = x,$

  \item[(b)] The process $\tilde{X}\eb := (X\ebx)_{x \in \RR^n}$ is a Feller Markov process on $\RR^n$ whose generator $\tilde{\cL}_{\beta}$ contains $\cC^2_c(\RR^n)$, the set of compactly supported $\cC^2$ functions, in its domain and such that for all $\tilde{f} \in  \cC^2_c(\RR^n)$,
  \bqn{tildeL}
 \nonumber \tilde{\cL}_{\beta}(\tilde{f}) &= &- \beta \langle \na \tilde{U}, \na \tilde{f} \rangle  -\f12\langle\na \wi U,\na\wi f \rangle+\f12\sum_{i=1}^n\wi E_i[\wi U]\wi E_i[\wi f]+ \tilde{U} \sum_{i = 1}^n \tilde{E}_i^2(\tilde{f})
  \\
  &=& - \beta \na\tilde{U}(\tilde{f}) -\f12\na\tilde{U}(\tilde{f}) +\f12\sum_{i=1}^n\wi E_i[\wi U]\wi E_i[\wi f]+ \tilde{U} \sum_{i = 1}^n \tilde{E}_i^2(\tilde{f}) \eqn
  \item[(c)] The map $x \mapsto X\ebx(t)$ is an homeomorphism. In particular,
$$\forall t \geq 0, \, X\ebx(t) \in \RR^n \setminus  \cU  \Leftrightarrow \exists t \geq 0, \,   X\ebx(t) \in  \RR^n \setminus \cU.$$
\end{description}
Set $S_i(x) =\sqrt{2\tilde{U}(x)}\tilde{E}_i(x).$ On $\RR^n\setminus \cU,$ (\ref{sde}) can be rewritten, using Stratonovich  formalism, as
\begin{eqnarray}
\label{sde-stra}
  dX(t) &=&  \left ( (-\beta - \f1{2}) \na \tilde{U}(X(t)) + \f1{2} \sum_{i = 1}^n DS_i(X(t)) S_i(X(t))  \right ) dt + \sum_{i = 1}^n S_i(X(t))  dB^{i}(t)\nonumber \\
   & = & (-\beta - \f1{2}) \na \tilde{U}(X(t)) + \sum_{i = 1}^n S_i(X(t)) \circ dB^{i}(t).
\end{eqnarray}
The vector fields $\nabla \tilde{U}$ and $S_i$'s being tangent to $N,$ this latter expression shows that $N$ (hence $M$) is invariant for $X\eb.$ That is:
$$\forall t \geq 0, \, X\ebx(t) \in N (\mbox{ resp. } M)  \Leftrightarrow \exists t \geq 0, \,   X\ebx(t) \in N (\mbox{ resp. } M).$$
It then follows that $X\eb:= (X\ebx)_{x \in M}$ is a Feller Markov process on $M,$ leaving $N$ invariant, whose generator $\cL_{\beta}$ contains $\cC^2(M)$ in its domain and such that $\cL_{\beta} f = \tilde{\cL_{\beta}} \tilde{f}|_M =  L_{\beta} f$ for all $f \in \cC^2(M)$ and $\tilde{f} \in C^2(\RR^n)$ such that $\tilde{f}|_M = f.$ The last equalities follows from  Lemma \ref{DeltaM} and (\ref{tildeL}), since on $M$ we have
\bq \langle\na \wi U,\na\wi f \rangle&=&\sum_{i=1}^n\wi E_i[\wi U]\wi E_i[\wi f]\eq

The strong Feller property on $N$ follows from the ellipticity of $L_{\beta}$ on $N$ (see e.g.\ Ichihara and  Kunita \cite{MR0381007,zbMATH03594399}, Lemma 5.1).
\subsection{On Remark \ref{rem4}}\label{5.2}
Given $0 < \lambda_- < \lambda_+,$ and $m \geq 2,$ let $D(\lambda_-,\lambda_+,m)$ be the set of diagonal matrices with entries
$\lambda_- = \lambda_1 \leq \lambda_2 \leq \ldots \leq \lambda_{m-1} \leq \lambda_m = \lambda_+.$ The set $\{\Lambda(A,\beta) \: : A \in D(\lambda_-,\lambda_+,m)\}$  is a compact interval $[\lambda_-(m,\beta), \lambda_+(m,\beta)]$ (as the image by  a continuous map of the compact connected set $D(\lambda_-,\lambda_+,m)$) contained in $[\lambda_-, \lambda_+].$

Let $A \in D(\lambda_-,\lambda_+,m)$ be the matrix with entries $\lambda_1 = \ldots \lambda_{m-1} = \lambda_-$ and $\lambda_m = \lambda_+.$
Then $$Z(\beta,A) = \int [\lambda_+ \theta_m^2 + \lambda_-(1-\theta_m^2)]^{-\beta} \sigma(d\theta)
= \EE \lt[\lt(\frac{\lambda_+X_m^2  + \lambda_- (\sum_{i = 1}^{m-1} X_i^2)}{\sum_{i = 1}^m X_i^2}\rt)^{-\beta}\rt],$$ where $X_1, \ldots, X_m$ are i.i.d. ${\cal N}(0,1)$ random variables.
By the strong law of large numbers and dominated convergence, this quantity converges, as $m \rightarrow \infty,$ toward $\lambda_-^{-\beta}.$
Thus $\lim_{m \rightarrow \infty} \lambda_-(m,\beta) = \lambda_-.$ Similarly, $\lim_{m \rightarrow \infty} \lambda_+(m,\beta) = \lambda_+.$

\subsection{On spherical integrals}\label{5.3}
In \eqref{Vlnr} we could have considered another  function $\sV$. Indeed, our first choice was
\bq
\wi\sV&\df& -\ln(U)\eq
since it seemed somewhat more ``intrinsic'' with respect to $U$. It can be shown similarly that the points [a] and [b] following \eqref{Vlnr} equally hold, with $\sV$ replaced by $\wi\sV$ and
$\sH_\beta$ by $\wi\sH_\beta$ given on $\{0\}\times\SS^{m-1}$ by
 \bq
 \fo \theta\in \SS^{m-1}, \qquad \wi\sH_{\beta}(0,\theta)&\df& -\trace(A)+2(1+\beta)\f{\lan \theta,A^2\theta\ran}{\lan \theta,A\theta\ran}\eq
where we recall that $A\df \Hess U(0)$.

\par
he sign of the quantity $\mu_{A,\beta}[\wi\sH_\beta(0,\cdot)]$ can then be used to discriminate between the attractiveness and repulsivity of $0$.
In particular  $\mu_{A,\beta}[\sH_\beta(0,\cdot)]$ and $\mu_{A,\beta}[\wi\sH_\beta(0,\cdot)]$ must have the same sign.
We tried to prove directly (without success!) that
\bqn{cur1}
 2(1+\beta)\mu_{A,\beta}[\phi_A]>\trace(A)&\Leftrightarrow& \beta>\f{m}2-1\\
\label{cur2}  2(1+\beta)\mu_{A,\beta}[\phi_A]<\trace(A)&\Leftrightarrow& \beta<\f{m}2-1\eqn
with
\bq
\fo \theta\in \SS^{m-1},\qquad \phi_A(\theta)&\df&  \f{\lan \theta,A^2\theta\ran}{\lan \theta,A\theta\ran}\eq

A by-product of our computations is thus to show the validity of \eqref{cur1} and \eqref{cur2}, which look as natural bounds on the corresponding spherical integrals for any given definite positive matrix $A$.

\vskip2cm
\hskip70mm
\vbox{
\copy4
 \vskip5mm
 \copy5
}


\begin{thebibliography}{1}

\bibitem{MR3155209}
Dominique Bakry, Ivan Gentil, and Michel Ledoux.
\newblock {\em Analysis and geometry of {M}arkov diffusion operators}, volume
  348 of {\em Grundlehren der Mathematischen Wissenschaften [Fundamental
  Principles of Mathematical Sciences]}.
\newblock Springer, Cham, 2014.

\bibitem{Ben18}
Michel Bena\"{\i}m.
\newblock Stochastic Persistence.
\newblock {\em arXiv}, 1806.08450, 2018.

\bibitem{BH22}
Michel Bena\"{\i}m and Tobias Hurth.
\newblock {\em Markov Chains on Metric Spaces
A Short Course},  Universitext.
\newblock Springer, 2022.

\bibitem{MR3910006}
Michel Bena\"{\i}m and Edouard Strickler.
\newblock Random switching between vector fields having a common zero.
\newblock {\em Ann. Appl. Probab.}, 29(1):326--375, 2019.

\bibitem{Echeverria}
Pedro Echeverria E.
\newblock A criterion for invariant measures of Markov processes.
\newblock {\em Z. Wahrscheinlichkeitstheorie und Verw. Gebiete.}, 61:1--16, 1982.

\bibitem{EK}
Stewart N. Ethier and Thomas G. Kurtz.
\newblock{\em Markov Processes, Characterization and Convergence}
\newblock John Wiley and Sons, 1986.

\bibitem{MR2088027}
Sylvestre Gallot, Dominique Hulin, and Jacques Lafontaine.
\newblock {\em Riemannian geometry}.
\newblock Universitext. Springer-Verlag, Berlin, third edition, 2004.

\bibitem{MR0381007}
   Kanji  Ichihara, and  Hiroshi Kunita.
    \newblock {\em A classification of the second order degenerate elliptic
              operators and its probabilistic characterization},
 \newblock {\em Z. Wahrscheinlichkeitstheorie und Verw. Gebiete.},30:235--254, 1974.

\bibitem{zbMATH03594399}
Kanji Ichihara and Hiroshi Kunita.
\newblock Supplements and corrections to the paper: {A} classification of the
  second order degenerate elliptic operators and its probabilistic
  characterization.
\newblock {\em Z. Wahrscheinlichkeitstheor. Verw. Geb.}, 39:81--84, 1977.

\bibitem{MR0885138}
   Wolfgang Kliemann.
   \newblock {\em Recurrence and invariant measures for degenerate diffusions},
   \newblock {\em Ann. Probab.}, 15(2):690--707, 1987.

\bibitem{MR1070361}
Hiroshi Kunita.
\newblock {\em Stochastic flows and stochastic differential equations},    volume 24 of
Cambridge Studies in Advanced Mathematics,
\newblock Cambridge University Press, Cambridge, 1990

\bibitem{MR3497465}
Jean-Fran\c{c}ois Le Gall.
   \newblock {\em Brownian motion, martingales, and stochastic calculus},
Graduate Texts in Mathematics,
\newblock Springer, Cham, 2016.




\bibitem{miclo:hal-04094950}
Laurent Miclo.
\newblock {On the convergence of global-optimization fraudulent stochastic
  algorithms}.
\newblock Preprint available at \texttt{https://hal.science/hal-04094950},
  April 2023.

\bibitem{Pennec_geometric}
Xavier Pennec.
\newblock Probabilities and statistics on {Riemannian} manifolds : A geometric
  approach.
\newblock Technical Report RR- 5093, inria-00071490, INRIA, 2004.

\bibitem{Stroock}
Daniel W. Stroock
   \newblock {\em An Introduction to the Analysis of Paths on a Riemmanian Manifold},
Mathematical Surveys and Monographs, Vol 74
\newblock American Mathematical Society, 2000.
\bibitem{enwiki:1157164633}
{Wikipedia contributors}.
\newblock List of formulas in {Riemannian} geometry --- {Wikipedia}{,} the free
  encyclopedia.
\newblock
  \url{https://en.wikipedia.org/w/index.php?title=List_of_formulas_in_Riemannian_geometry&oldid=1157164633},
  2023.
\newblock [Online; accessed 15-October-2023].


\end{thebibliography}
\end{document}